\documentclass{article}
\usepackage[utf8]{inputenc}

\title{Permanence as a Principle of Practice}

\author{Iulian D. Toader\thanks{\bigskip Institute Vienna Circle, University of Vienna, iulian.danut.toader@univie.ac.at}}

\date{}

\begin{document}

\maketitle

\begin{flushright}
Dedicated to the memory of Mic Detlefsen
\end{flushright} 

\bigskip \bigskip \bigskip

The paper discusses Peano's argument for preserving familiar notations. The argument reinforces the principle of permanence, articulated in the early 19th century by Peacock, then adjusted by Hankel and adopted by many others. Typically regarded as a principle of theoretical rationality, permanence was understood by Peano, following Mach, and against Schubert, as a principle of practical rationality. The paper considers how permanence, thus understood, was used in justifying Burali-Forti and Marcolongo's notation for vectorial calculus, and in rejecting Frege's logical notation, and closes by considering Hahn's revival of Peano's argument against Pringsheim' reading of permanence as a logically necessary principle.

\bigskip

Der Beitrag diskutiert Peanos Argument für die Bewahrung vertrauter Notationen. Es geht um das Prinzip der Permanenz, das im frühen 19. Jahrhundert von Peacock artikuliert, dann von Hankel angepasst und von vielen anderen übernommen wurde. Während es im allgemeinen als ein Prinzip der theoretischen Rationalität betrachtet wurde, interpretierte Peano es anknüpfend an Mach und gerichtet gegen Schubert als ein Prinzip der praktischen Rationalität. Der Beitrag untersucht, wie Permanenz, so verstanden, von Burali-Forti und Marcolongo zur Rechtfertigung von Bezeichnungen in der Vektorrechnung verwendet wurde und zur Zurückweisung der logischen Notation Freges führte. Die Arbeit schließt mit einer Erörterung der Wiederbelebung des Peanoschen Arguments in Hahns Kritik von Pringsheims Interpretation der Permanenz als ein logisch notwendiges Prinzip.

\bigskip


\section{Introduction}

Should familiar notation be preserved in science? Should one continue to use the symbols used in already established theories? And should one commit to preserving the properties these symbols possess in such theories? The history of science does not seem to provide a definite answer, as indicated on the one hand by the existence of alternative notations for one and the same theory, and on the other hand by unequivocal attempts to increase the uniformity of notation across theories. In this paper, I discuss an argument in favor of the preservation of familiar notation in logic and mathematics, explicitly advocated as such, early in the 20th century, by Giuseppe Peano. 

Peano's argument reinforces a rather different understanding of the principle of permanence than that articulated in the first part of the 19th century by George Peacock and then adjusted and advocated by Hermann Hankel, and subsequently adopted by many other mathematicians and philosophers of mathematics. While they took it to be a principle of theoretical rationality, i.e., one that, albeit not universally valid, is nevertheless indispensable for the development of mathematics as a genuine science, the principle of permanence was understood by Peano, following Mach, as a principle of practical rationality. More specifically, it was understood to justify the preservation of familiar notation in science, as illustrated by Cesare Burali-Forti and Roberto Marcolongo's notational project for vectorial calculus. By contrast, Peano regarded Frege's innovative notation for logic as an obstruction or an injunction against permanence thus understood. As will become clear, one salient difference between these two conceptions of the principle of permanence is that, while as a principle of theoretical rationality it was applied to the basic laws and rules of a mathematical theory, as a principle of practical rationality it was applied more generally to ideas and relations (of which the basic laws and rules of mathematical theories are but a subset), as well as to their notational expressions. Another related difference concerns the goals motivating the application of the principle of permanence: preservation of basic laws and rules was motivated by the goal of rendering an extension of a mathematical theory as meaningful as the theory undergoing the extension, while preservation of familiar ideas and notations was motivated by the ideal of thought economy. 

Is Peano's conception of the principle of permanence justified? What exactly did he find problematic with its understanding as a principle of theoretical rationality? What made him adopt an alternative conception, and why exactly did he associate this with the Italian project for vectorial calculus? What made Peano think that Frege's logical diagrams interfere with the principle of permanence, specifically understood as a principle of practical rationality? Finally, is the principle of permanence, thus understood, sufficient to eliminate all arbitrariness in the development of new theories?

The paper will proceed as follows. In section two, I recall some of the relevant historical background to the principle of permanence, spanning almost a century, from the 1820s to the 1910s, including an exchange between Augustus De Morgan and Ada Lovelace, so far underexplored. This will provide evidence that not everybody understood this principle in the same way: some thought it was a methodological principle, some others a metaphysical one; some took it to represent an indispensable foundation for science, some others a necessary truth that could be assumed in mathematical proofs. Yet others, like Mach, considered the principle of  permanence a matter of practical rationality, meant (just like logical consistency) to improve thought economy and reduce intellectual discomfort. In section three, after briefly describing the view of permanence defended by Hermann Schubert, along with the criticism that Frege stormed against it, I reconstruct Peano's argument, also directed against Schubert, for the claim that permanence could be regarded only as a principle of practical rationality. Then I consider what Peano took to be an illustration of the principle, thus understood, in the notational project initiated by Burali-Forti and Marcolongo, where it was specifically meant to facilitate the dissemination and assimilation of vectorial calculus. 

Peano's conception of permanence implies that the introduction of a new notation in science, which fails to preserve relevant properties of familiar notations already in use, would be costly from a practical point of view. Of course, what would be the relevant properties that ought to be preserved, and the practical costs for failing to do so, is to be determined in each particular case. But for vectorial calculus, which will be here discussed, the claim was that the larger the variety of notations in use, the more difficult the dissemination and the assimilation of the calculus. For another case in point, I discuss in section four Peano's earlier rejection of the two-dimensional notation that Frege had introduced in his \textit{Begriffsschrift}. As is well known, Frege's sweeping criticism of formalism repudiated what he took to be Hankel's formalism and, together with this, the principle of permanence as an indispensable foundation for formal mathematics. But Frege's dispute over notation with Peano (and Schr\"oder) indicates that he rejected the principle of permanence as a practical concern as well. It also indicates that at least some of their criticisms of the \textit{Begriffsschrift} were not simply based on misunderstanding its innovative character and its higher standards of rigor, but were also motivated (perhaps more explicitly in Peano's case) by a conception of permanence as a principle of practical rationality.

Despite their opposite character relative to the principle of permanence, both Frege's and Burali-Forti and Marcolongo's notational projects failed. Moreover, Schubert's understanding of permanence as a principle of theoretical rationality not only survived Peano's criticism, but in fact received purportedly stronger support, such as in the works of Alfred Pringsheim. This led in turn to a reiteration of Peano's criticism, by Hans Hahn, which I discuss in section five. Against Pringsheim's insistence on permanence as a principle that allows no choice as to the rules to be preserved when extending a number system, Hahn argued that even if permanence were a logically necessary principle, arbitrariness would still creep in, since the principle of permanence cannot rule out the possibility of distinct extensions of a number system. In order to eliminate this arbitrariness, he maintained that a further condition is required to guarantee extension uniqueness. As we will see, Hahn goes beyond Peano by arguing that in order to solve this problem of arbitrariness in the extension of number systems, what is needed, beside permanence, is the stipulation that equality must be uniquely defined in the extended system, which in turn requires that its basic operations be conceived as one-to-one mappings. He also noted, however, that practical considerations might also help sometimes to choose between distinct extensions.

\section{The  Principle of Permanence}

Implicit mention of what would come to be known as the principle of permanence can be found early in the 19th Century, for example in Johann Peter Wilhelm Stein's 1828 \textit{Elemente der Algebra}: ``One of the most important difficulties that calculating with numbers has to overcome is that it is to establish rules for those new cases in calculation which are \textit{as much as possible} in agreement with the customary arithmetical rules and still properly well founded. [...] Of the rules of arithmetic, which have indeed been proven for all possible number values of letters, but also only for number values, one has to wish that they were also right if the letters occurring in them, all or partially, meant zero or negative quantities.''\footnote{English translation quoted from Schubring (2005, 514, my emphasis). For discussion of the principle of permanence, see e.g. Koppelman (1971), Pycior (1997), Fisch (1999), and Detlefsen (2005). As a general rule for this paper, I give my own English translations in the text only if they are not already available, and quote the original versions in footnotes.} Thus, Stein thought that although one must preserve the rules of the arithmetic of positive natural numbers when extending the range of computations to include, e.g., negative natural numbers, the new rules for an extended system of numbers must be in agreement with the old ones ``as much as possible''. This suggests that, on Stein's view, the new rules must be either identical to the old ones, or, if different, must be at least consistent with them, and raises a question as to how much disagreement with the old rules may be allowed such that one may still call the new rules ``properly well founded.'' As we will see in section three, a similar question would come up in relation to later expressions of the principle of permanence, critically discussed by Peano.

An explicit formulation of the principle of permanence appears in George Peacock's 1833 ``Report on the recent progress and present state of certain branches of analysis'', where it is expressed in the following way: ``Whatever form is algebraically equivalent to another when expressed in general symbols, \textit{must continue to be} equivalent, whatever those symbols denote.'' This ``Principle of Permanence of Equivalent Forms'', as Peacock called it, also stipulates the following: ``Whatever equivalent form is discoverable in arithmetical algebra considered as the science of suggestion, when the symbols are general in their form, though specific in their value, \textit{will continue to be} an equivalent form when the symbols are general in their nature as well as in their form.''\footnote{See Peacock (1833, 198sq, my emphasis). These same formulations are repeated in Peacock's two-volume \textit{Treatise on Algebra} (Peacock 1842-1845).} For example, if \textit{m, n} and \textit{a} denote any integers, as they do in what Peacock called arithmetical algebra, then it is the case that $ma+na=(m+n)a$. Arithmetically equivalent forms like this one are said to have a ``necessary'' existence (Peacock 1833, 199), for they are true in virtue of the definitions of the basic operations (like addition and multiplication). Permanence demands that $ma+na=(m+n)a$ remains an equivalent form in what Peacock called symbolic algebra, where \textit{m, n} and \textit{a} may denote anything whatsoever. However, as a symbolically equivalent form, $ma+na=(m+n)a$ cannot be true in virtue of the definitions of the operations, since in symbolic algebra no such definitions are yet available. The meaning of the basic operations is only determined by symbolically equivalent forms like $ma+na=(m+n)a$. Such forms are thus said to have a ``conventional'' existence (Peacock 1833, 200).

One of Peacock's main concerns was the applicability of symbolic algebra. He noted that although the rules of symbolic algebra, i.e., a certain set of symbolically equivalent forms, can be used for deducing other symbolically equivalent forms, i.e., the theorems of symbolic algebra, this would be a mere game with meaningless symbols if the forms had no applications, that is, if they had \textit{only} a conventional existence. Thus, Peacock required that symbolically equivalent forms permit an arithmetical interpretation, i.e., that they can always be transferred back to arithmetical algebra, such that they can also have a necessary existence as arithmetically equivalent forms. Another of Peacock's main concerns was the generality of the principle or permanence. As stated above, the principle requires that all arithmetically equivalent forms should be preserved in symbolic algebra, and vice versa, all symbolically equivalent forms should be transferable back to arithmetical algebra. But Peacock knew very well that some arithmetically equivalent forms are essentially connected to the specific values of symbols, such that they cannot be preserved in symbolic algebra. Some arithmetically equivalent forms are, therefore, only hypothetically preserved in symbolic algebra, in the sense that as symbolically equivalent forms they have only a ``hypothetical'' existence (Peacock 1833, 210): they always ``degenerate'' into arithmetically equivalent forms that hold for some specific numerical values of their symbols. Due to such cases, Peacock recommended caution in the use of the principle of permanence. Despite such limitations of its generality, he conceived of it as a principle of theoretical rationality. That is, he thought that permanence, just like logical consistency, was indispensable to the development of symbolic algebra as a genuine science: the idea that equivalent forms are to be preserved ``must guide'' this development. Thus, albeit not universally valid, permanence was taken to be a methodologically necessary principle. Peacock implied that, without taking permanence as our guide, we would end up with a set of arbitrary equivalent forms, with no application and no meaning whatsoever. That, he thought, would hardly be deserving of the name of science.

An arguably similar conception of permanence was later defended by Hermann Hankel in his 1867 book, \textit{Vorlesungen \"uber die complexen Zahlen und ihre Functionen}. His development of purely formal theories of numbers, disconnected from intuition and constrained only by the conditions of logical consistency and mutual independence of its rules, was to be similarly guarded against potential meaninglessness. To this effect, Hankel required that the formal rules for the operations with what he called objects of thought (\textit{Gedankendinge}) admit as their subordinate the non-formal rules for the operations with quantities as objects of intuition, e.g., the rules of arithmetic. This requirement was meant to ensure that the statements of his formal theories would always have a non-formal, arithmetical interpretation and, thus, applicability. Hankel wrote: ``Any attempt to treat the irrational numbers formally and without the concept of quantity must lead to extremely abstruse and laborious artificialities that, even if they can be carried out with complete rigor, which we have good reason to doubt, do not have a higher scientific value. Because in general it is a matter of systematic science to get clear on and be aware of the true foundations of the natural development of ideas, but not to want to replace the organism with its always fresh production power by a lifeless and unproductive mechanism, even if ingeniously constructed.''\footnote{``Jeder Versuch, die irrationalen Zahlen formal und ohne den Begriff der Gr\"osse zu behandeln, muss auf h\"ochst abstruse und beschwerliche K\"unsteleien f\"uhren, die, selbst wenn sie sich in vollkommener Strenge durchf\"uhren liessen, wie wir gerechten Grund haben zu bezweifeln, einen h\"oheren wissenschftlichen Wert nicht haben. Denn überall ist es Sache der systematischen 
Wissenschaft, sich der wahren Grundlagen der natürlichen Entwickelung der Ideen klar und bewusst zu werden, nicht aber den Organismus mit seiner immer frischen Productionskraft durch einen, wenn auch scharfsinnig construirten, doch todten und unproductiven Mechanismus ersetzen zu wollen.'' (Hankel 1867, 46sq)} Hankel would later be criticized for imposing such a requirement on formal theories. Alfred Pringsheim, for example, to whom I shall return below in section five, commented on the passage just quoted: ``It appears extremely remarkable that even the creator of a \textit{purely formal theory of rational numbers} has shown so little understanding for the corresponding further development of the number concept.''\footnote{``Es erscheint \"ausserst merkw\"urdig, dass gerade der Sch\"opfer einer \textit{rein formalen Theorie der Rationalzahlen} f\"ur die entsprechende Weiterbildung des Zahlbegriffs so wenig Verst\"andnis gezeigt hat.'' (Pringsheim 1894, 57, n32) To be clear, the criticism here is not that Hankel took the principle of permanence to be universally valid, that is, that it did not admit of any exceptions; rather, Pringsheim's criticism is that a purely formal theory should not be constrained by any applicability requirements.} 

Hankel gave the principle of permanence the following formulation: ``If two forms expressed in general signs of the universal arithmetic are equal to one another, they should remain equal if the signs cease to denote simple quantities and the operations thereby receive some different content as well.''\footnote{``Wenn zwei in allgemeinen Zeichen der arithmetica universalis ausgedrückte Formen einander gleich sind, so sollen sie einander auch gleich bleiben, wenn die Zeichen aufhören, einfache Grössen zu bezeichnen, und daher auch die Operationen einen irgend welchen anderen Inhalt bekommen.'' (Hankel 1867, 11)} This corresponds roughly to Peacock's formulation of the principle: $ma+na=(m+n)a$, as an arithmetical expression of equal forms should be preserved in formal mathematics, where \textit{m, n} and \textit{a} may denote any objects of thought and $ma+na=(m+n)a$ is taken as a formal expression. Hankel thought that permanence was not a mere guide or a merely heuristic (\textit{hodegetische}) principle; rather, he claimed that it was indispensable, for he thought that it was a ``metaphysical'' principle (1867, 12). However, Hankel also clearly warned against an incautious universal application of the principle, pointing out that in developing formal theories, certain rules or laws that hold for the real numbers cannot be extended to complex and hypercomplex numbers. Furthermore, he proved that there can exist no  extension beyond the complex numbers that preserves the commutativity of basic operations (Detlefsen 2005, 286). 

Peacock and Hankel converged, I think, on the point that the principle of permanence, although not universally valid, is to be understood as an indispensable guide for the development of symbolic algebra and of formal number theories, that is, as a principle of theoretical rationality, primarily meant to guarantee the arithmetical interpretability of such theories, and thus to justify the claim that they are genuine sciences in their own right, rather than meaningless games.\footnote{For a different reading of Hankel relative to Peacock, see Peckhaus (1997).} But other mathematicians held stronger views about permanence. 

Augustus De Morgan, for instance, took the principle of permanence to be not merely an indispensable guide, but a necessary mathematical truth. His own formulation of the principle is as follows: ``all algebraical expressions are combined and reduced by rules, which, although derived from notions on quantity, will produce the same results, if we alter the form of the primitive expressions in any manner, consistently with the rules, even though the new forms should no longer admit of being considered as quantities.''\footnote{De Morgan (1836-1842, 119), quoted in Hollings \textit{et al.} (2017, 18). De Morgan also formulated the principle of permanence as follows: ``When an algebraical multiplication, or other operation, such as has hitherto been defined, can be proved to produce a certain result in cases where the letters stand for whole numbers, then the same result \textit{must be true} when the letters stand for fractions, or incommensurable numbers, and also when they are negative.'' (1837, 212, quoted in Hollings \textit{et al.} 2017, 17)} What De Morgan seems to say is that the forms of algebraic expressions, which result from the application of a set of rules to primitive expressions, remain the same, whether the forms admit of being considered as quantities or not, i.e., no matter what the expressions denote. He seems to imply that all rules applicable to expressions the forms of which admit of being taken as quantities are preserved when the forms to be altered belong to expressions that do not denote quantities any longer, and so that the principle would be generally applicable. 

In correspondence with Ada Lovelace, who had questioned the validity of this principle, De Morgan insisted on ``the necessity of its truth.''\footnote{Quoted in Hollings \textit{et al.} (2017, 18). De Morgan's full reaction to Lovelace's skepticism was the following: ``This principle requires some algebraical practice to see the necessity of its truth.'' (The Lovelace Byron Papers, Bodleian Library, Oxford, Box 170, Transcripts of Folios 1-179, by Christopher Hollings. Available on the Clay Mathematics Institute's website: https://claymath.org)} His use of the principle of permanence in mathematical proofs, such as his proof of the binomial theorem, suggests that he attributed to the principle a necessity stronger than methodological necessity. Unconvinced, Lovelace followed up: ``It cannot help striking me that this extension of Algebra ought to lead to a further extension similar in nature, to Geometry in Three-Dimensions; \& that again perhaps to a further extension into some unknown region, \& so on ad-infinitum possibly.'' (quoted in Hollings \textit{et al.} 2017, 18) It's not quite clear what Lovelace suggested here, and whether she thought that a potentially infinite succession of extensions ``similar in nature'' was problematic, and if so, why. But one can speculate that she might have meant that, if the principle of permanence is considered a necessary and universally valid truth, then all possible extensions of a given mathematical system would be essentially the same. This might have been what she found ``striking''. If this is what she meant, then she anticipated an intuition that, as we will see in the next section, played an important part in Peano's argument against permanence as a principle of theoretical rationality. 

In any case, De Morgan's understanding of the principle of permanence was much closer to Peacock's than was the view later advanced by Mach, who saw it as a matter of practical, rather than theoretical rationality. In his 1905 book \textit{Knowledge and Error}, he described a view of permanence as yet another handmaid of thought economy: ``The mutual adaptation of thoughts is not exhausted in the removal of contradictions: whatever divides attention or burdens the memory by excessive variety, is felt as uncomfortable, even when there are no contradictions left. The mind feels relieved whenever the new and unknown is recognized as a combination of the known, or the seemingly different is revealed as the same, or the number of sufficient leading ideas is reduced and they are arranged according to the principles of permanence and sufficient differentiation. Economizing, harmonizing and organizing of thoughts are felt as a biological need far beyond the demand for logical consistency.'' (Mach 1905, 127sq) Mach understood the principle of permanence as a principle of practical rationality, subordinate to his general principle of thought economy, in the same manner in which he conceived of consistency: together, they were meant to achieve the supreme goal of mutual adaptation of thoughts (including the adaptation of thoughts to facts). The special job of permanence was to eliminate, as much as possible, the intellectual discomfort caused by loading our memory and dividing our attention through new ideas. But although Mach may have taken it to be applicable in all generality to ideas, the principle or permanence was more specifically to be applied to ``leading ideas''. Commenting on Mach's view of laws as the leading ideas that  economically order our experiences, Musil took such ideas to ``correspond to the need for permanence'' and justified this claim in the following way: ``For it is in them -- in constant laws and equations [...] that thought seeks to grasp those ideas which can be held on to permanently whatever individual changes may occur.'' (Musil 1908, 24) Thus, the principle of permanence, although it could be generally applied to ideas, was to be applied more particularly to laws and equations, and as such it was supposed to eliminate (or reduce) intellectual discomfort not only by preserving old and familiar laws, but also by seeing that new laws, which could not be seen as combinations of the old ones, are avoided as much as possible.\footnote{For a discussion of Mach's conception of the principle of permanence as a thought-economical principle, and his dispute with Husserl on the matter, see Toader (2019).}

To better understand Mach's view of permanence, it may be helpful to consider it as against a Peircean conception of the nature of inquiry. As is well known, Peirce defended a view according to which inquiry, and in particular scientific inquiry, lacks proper motivation in the absence of logical inconsistency. He famously wrote: ``That the settlement of opinion is the sole end of inquiry is a very important proposition. It sweeps away, at once, various vague and erroneous conceptions of proof. Some people seem to love to argue a point after all the world is fully convinced of it. But no further advance can be made. When doubt ceases, mental action on the subject comes to an end; and, if it did go on, it would be without a purpose.'' (Peirce 1877, 11) Purposeful mental action and, thus, proper scientific inquiry requires the presence of ``real and living'' (as opposed to mere Cartesian) doubt, which arises only when logical inconsistencies are revealed. This doubt stimulates the mind to seek the settlement of belief. Once inconsistency is removed, stimulation ceases, belief settles, and inquiry should stop. To postpone settlement, which is what Peirce thought was characteristic of Cartesian epistemology, is a ``perversity.'' (Peirce 1878, 39) But to seek improvement after belief has settled would be no less irrational. Still, this is precisely what Mach seems to think one should seek: as we have just seen, he thought that getting rid of logical inconsistencies, and thus reaching a Peircean settlement of belief, would never be the end of any inquiry, as it may not be enough for reaching intellectual comfort. What Mach further required is that old ideas, laws and equations, be preserved as much as possible, and that new ideas, laws and equations, irreducible to the old ones, be avoided as much as possible. Only if this extra two-fold requirement were satisfied, could one fully achieve a mutual adaptation of thoughts.

\section{Peano's Argument}

Before we analyze Peano's argument about the principle of permanence, we should briefly discuss Schubert's conception of permanence as a principle of theoretical rationality, for this is what Peano reacted to. Schubert contributed to the very first volume of the \textit{Encyklop\"adie der mathematischen Wissenschaften}, a joint editorial project of the academies of sciences in G\"ottingen, Leipzig, Munich, and Vienna. Citing both Peacock and Hankel, Schubert invoked the principle of permanence or, as he also sometimes called it, the principle of exceptionlessness (\textit{Ausnahmslosigkeit}), each time a new extension of the system of positive natural numbers was introduced. His formulation of it has four parts: ``1. to give to each concatenation of signs, which represents no previously defined number, such a sense that the concatenation can be manipulated after the same rules as if it represented one of the previously defined numbers; 2. to define such a concatenation as a number in an extended sense of the word, and thereby to extend the concept of number; 3. to prove that for the numbers in the extended sense, the same propositions hold as for the numbers in the not yet extended sense; 4. to define what equal, greater than, and less than, mean in the extended number domain.''\footnote{``[D]as \textit{Prinzip der Permanenz} [...] in viererlei besteht: erstens darin, jeder Zeichen-Verkn\"upfung, die keine der bis dahin definierten Zahlen darstellt, einen solchen Sinn zu erteilen, dass die Verkn\"upfung nach denselben Regeln behandelt werden darf, als stellte sie eine der bis dahin definierten Zahlen dar; zweitens darin, eine solche Verkn\"upfung als Zahl im erweiterten Sinne des Wortes zu definieren und dadurch den Begriff der Zahl zu erweitern; drittens darin, zu beweisen, dass f\"ur die Zahlen im erweiterten Sinne dieselben S\"atze gelten, wie f\"ur die Zahlen im noch nicht erweiterten Sinne; viertens darin, zu definieren, was im erweiterten Zahlengebiet gleich, gr\"osser und kleiner heisst.'' (Schubert 1898, 11)} Although the relations between these four parts are never clarified, the same formulation is given again in Schubert's 1899 book, \textit{Elementare Arithmetik und Algebra}. What is quite clear, however, at least from the first and the third part of the principle, is that he thought that all rules that govern a non-extended number system must govern the extended system as well, and furthermore, that all propositions that hold in the non-extended system must provably hold in the extended one as well. 

Not everybody agreed that this was clear enough, though. Schubert's contribution to the \textit{Encyklop\"adie} was reviewed by Frege, whose discussion of it is a prime example of the kind of contempt and derision that he was capable of spewing at some of his colleagues (Frege 1899). It may also be considered a prime example of the kind of hurried and unfair reading that Frege gave to some texts. For example, he read into Schubert's first part of the principle the idea that the rules follow from the sense of the signs, although Schubert does not actually say that; rather, what he said is that the application of the rules requires that the signs have a sense. Frege then argued that to believe that the propositions which hold for the new numbers are identical to the propositions that hold for the old ones \textit{because} the rules for the manipulations of the new number-signs are identical to the rules for the manipulations of the old number-signs proves that Schubert confused numbers with their signs, which on Frege's view was of course a fatal sin. However, as already mentioned, Schubert unfortunately never clarified the relations between the four parts of his statement of the principle of permanence, and in particular between the first and the third one, so it is not clear that he actually believed what was attributed to him. Thus, it's not clear that Schubert was as much of a sinner as Frege took him to be. More importantly, I think that this is indicative of Frege's attitude towards permanence as a principle of theoretical rationality. In the next section, further below, after briefly recalling his evisceration of Hankel, I will discuss what Frege seems to have thought of permanence, as Peano understood it, i.e., as a principle of practical rationality.

Unlike Frege, and a decade after, Peano directed his criticism of the principle of permanence solely at the third part of Schubert's formulation:\footnote{Peano translated Schubert's full formulation of the principle in an earlier paper, written in his simplified Latin (Peano 1903).}

\begin{quote}
    
This principle of permanence reached its apogee with Schubert, who, in the \textit{Encyklop\"adie der mathematischen Wissenschaften}, affirmed that one must ``prove that for the numbers in the extended sense, the same propositions hold as for the numbers in the not yet extended sense.'' Now, if all the propositions which are valid for the entities of one category are valid also for those of a second, then the two categories are identical. Hence --- if this could be proved --- the fractional numbers are integers! In the French edition of the \textit{Encyklop\'edie} these things are put to rights. There it says that one must be ``guided by a concern for keeping the formal laws as much as possible.'' Thus the principle of permanence acquires the value of a principle, not of logic, but of practice, and is of the greatest importance in the selection of notation. Basing their work on precisely this principle --- a particular case of what Mach called the principle of economy of thought --- Professors Burali-Forti and Marcolongo succeeded in untangling the disordered skein of notations in vectorial calculus, where all used to be arbitrary (and many still believe that the notations are necessarily arbitrary). (Peano 1910, 225) 

\end{quote}

The argument presented here by Peano could be initially reconstructed as follows: (P1) It is not the case that all propositions that are true in an extended system are provably true in the non-extended one. (P2) Hence, only the formal laws of the non-extended system must be preserved as much as possible in the extended one.  (P3) Thus, the principle of permanence is not a logical one, but only a principle of practice. (P4) Burali-Forti and Marcolongo interpreted it in precisely this way in their notational project for the vectorial calculus.

I will briefly consider the project mentioned in (P4) further below. For now, however, let us try to clarify Peano's argument. The first premise (P1) appears to have been justified by the following \textit{reductio}: if all propositions that hold in an extended system provably hold in the non-extended system as well, then the two number systems would be identical, which would be a contradiction (for it would obviously deny the assumption that an extension is different than the non-extended system). But even if one accepted this justification for (P1), it's not clear how (P2) would follow; it's not clear, that is, how it would follow that the formal laws of the non-extended system must be preserved as much as possible, rather than completely, in the extended system (Detlefsen 2005, 287). For it is of course possible that although some propositions that hold in the extended system are not true in the non-extended one, the formal laws of the non-extended system are all preserved. This, however, might not be how Peano intended his argument. For it's not clear that he intended (P2) to follow from (P1), that he intended the qualified version of the principle to follow from the rejection of the unqualified one. Indeed, he might have considered (P1) and (P2) as strictly unconnected, corresponding to distinct stages in Schubert's own formulation of permanence. If so, then Peano's argument could be reconstructed as follows: (P1') The formal laws of the non-extended system must be preserved as much as possible. (P2') Thus, the principle of permanence is not a logical one, but only a principle of practice. (P3')  Burali-Forti and Marcolongo interpreted it in precisely this way in their notational project for the vectorial calculus.

Drawing on Peano's reading of the French translation of Schubert's contribution, the justification for (P1') most plausibly stems from views on permanence such as those we have already seen articulated by Stein, even before Peacock and Hankel. Such views were sensitive to the fact that, in some cases, one has no choice but to forgo some formal laws when extending a number system. Inferring from this that permanence cannot be a logical principle, as stated by (P2'), Peano clearly meant that (P1') was sufficient to reject views like De Morgan's, who understood the principle of permanence as a necessary truth, which can be justifiably used in proofs. But did Peano also mean to reject Peacock's and Hankel's views about it, as an indispensable guide for the development of mathematics? I think that he did, to the extent that there is no precise account of what ``as much as possible'' means in (P1'). For this seems to allow that the choice as to what formal laws to preserve, and what to forgo, is logically arbitrary, in the sense that more than one option is logically possible. But then one can, at least in principle, choose not to preserve any formal laws whatsoever, and postulate instead entirely new ones. This renders permanence dispensable, insofar as one could develop, say, a new theory of numbers independently of any available theories of numbers. But instead of dismissing such a theory as meaningless and unscientific, like Peacock and Hankel would have done, Peano thought that it should be rejected because its development would have a high practical cost: its dissemination and assimilation would require a tremendous intellectual effort. He also thought that embracing the principle of permanence could help avoid this cost. This is why he adopted a Machian conception of permanence, to replace what he took to be a defective conception of it as a principle of theoretical rationality. 

Thus understood by Peano, permanence was considered useful, for example, in justifying various attempts to make notations more uniform across science. He expressed this idea in an early version of his famous 1921 paper on definitions in mathematics: ``Another practical law that governs mathematical definitions and notations is the so-called principle of permanence, or principle of conservation of the formal laws. This principle requires that when establishing a new system of notations, or a new calculus, it is convenient to do it such that the new calculus be similar \textit{as much as possible} with old calculi, so that the student does not have to learn a whole new calculus, but only the differences from the theory known to him. This principle corresponds to the principle of the economy of thought, and of minimum work.''\footnote{``Un altra legge pratica che regge le definizioni e le notazioni matematiche, è il cosi detto principio di permanenza, o principio di conservazione delle leggi formali. Questo principio impone che quando si stabilisce un nuovo sistema di notazioni, o un nuovo calcolo, conviene fare in modo che il nuovo calcolo sia simile \textit{per quanto è possibile} ai calcoli antichi, sicchè lo studioso non debba imparare tutto un nuovo calcolo, ma solo le differenze colla teoria a lui nota. Questo principio corrisponde al principio dell'economia del pensiero, e del minimo lavoro.'' (Peano 1911, 69, my emphasis) The English translation of Peano (1921), included in Peano (1973), does not contain this passage.} Peano thought that this conception of permanence was suitably illustrated by Burali-Forti and Marcolongo's notational project for vectorial calculus. We will briefly describe this project below, and then explain why Peano also believed that Frege's logical notation went against the principle of permanence.

\section{Permanence and Notation}

On Peano's conception of the principle of permanence, the existence of a plurality of notations for the same theory, as well as the proposal of a new notation that fail to preserve some (or any) characteristics of familiar notations already in use, would be regarded as a major drawback from a practical point of view. Peano thought that one should want to avoid this, and as we will presently see, this is at least part of the reason why he rejected Frege's notation for logic, and also part of the reason why he strongly supported Burali-Forti and Marcolongo's notational project for vectorial calculus. Let us start with the latter.

In the first decade of the 20th century, there existed many different notations in use for even the basic or fundamental notions of the vectorial calculus, sometimes within the works of the same mathematician. Among Italian mathematicians, in particular, this state of affairs was perceived as ``anarchical'' (Sallent Del Colombo 2010, 514), which motivated the need to find a ``unique and universal'' notation. In their 1909 book, Burali-Forti and Marcolongo deplored the fact that Hamilton used different combinations of symbols for denoting vectors, e.g., ``B--A'' and then ``AB'' for a vector from point A to point B. They further noted that others used different symbols for vectors: ``The AB notation does not agree with the simple notation used by Grassmann for the geometric formations [...] of which the vectors are a special case.'' But more importantly, they continued, Hamilton's notation ``destroys the analogy with algebraic calculus, that is, it is contrary to the useful principle of permanence."\footnote{``La notazione AB non è d'accordo con la notazione semplice usata da Grassmann per le formazioni geometriche [...] delle quali i vettori sono un caso particolare. Distrugge la analogia con il calcolo algebrico, cioè è contraria all'utile princìpio di permanènza.'' (Burali-Forti and Marcolongo 1909, 240)} Burali-Forti and Marcolongo considered the analogy between the vectorial and the algebraic calculus as a consequence, or a particular expression, of the principle of permanence. They also thought that a notation for vectorial calculus that was meant to be unique and universal must be in agreement, as far as possible, with the notations already developed for that calculus. 

The two requirements for what Burali-Forti and Marcolongo thought would be a successful development of their project may, therefore, be summarized as follows: ``1) the notations (at least the fundamental ones) must not be in contradiction with those (also fundamental ones) of Möbius, Hamilton, Grassmann [...]; 2) the vector algorithm must be established in such a way as to deviate \textit{as little as possible} from the universally known algorithm of algebra, because respecting the laws of permanence and economy greatly facilitates the dissemination of the vectorial calculus.''\footnote{``1) le notazioni (almeno le fondamentali) non devono essere in contraddizione con quelle (pure fondamentali) di Möbius, Hamilton, Grassmann, perché, anche avuto riguardo al sistema vettoriale minimo occorrente in pratica, non pare lecito ipotecare il passato e l'avvenire delle grandi opere di quei grandi;
2) l'algoritmo vettoriale deve essere stabilito in modo da discostarsi \textit{il meno possibile} da quello universalmente noto dell'algebra, perché rispettando le leggi di permanenza e di economia si facilita grandemente la diffusione del calcolo vettoriale.'' (Sallent Del Colombo 2010, 515sq, my emphasis)} These requirements are further articulated and justified in a series of articles published in \textit{Rendiconti del Circolo Matematico di Palermo}, where Burali-Forti and Marcolongo chose to present their project in detail. Upon the publication of the articles, a debate on their proposed notation ensued in \textit{L’Enseignement Mathématique}, which included among the contributors Felix Klein, who was rather skeptical that the project would succeed, and Peano, who defended it wholeheartedly. Nevertheless, the project ended up by failing to meet everyone's standards. The American mathematician Edwin Bidwell Wilson --- a former student of Gibbs at Yale --- concluded for example that ``there is no apparent gain in uniformity of notations attributable to these Italian activities.'' (Wilson 1913, 525)

Peano's reasons for defending Burali-Forti and Marcolongo's project, and for endorsing their use of the principle of permanence in justifying the project, brings to light some neglected aspects of his criticism of Frege's innovative notation for modern logic. These aspects concern the conspicuous differences between Frege's notation and the notations already in use by most logicians, in particular Boole and his followers, but also its perceived disanalogy with calculi developed elsewhere in mathematics, in particular in arithmetic. But before we consider Peano's and Frege's different attitudes towards permanence as a principle of practical rationality, we should note Frege's misgivings with Hankel's view of permanence as a principle of theoretical rationality. 

Frege's reading of Hankel, just like his review of Schubert's \textit{Encyklop\"adie} contribution, seems to have been unfortunately rather inaccurate, but this did not stop him from dismissing most of what Hankel did as a mathematician, and treating it like ``a rhetorically useful attack-target-quote generating machine.'' (Tappenden 2019, 238sq) For example, just like he would later do against Schubert, Frege blamed Hankel for allegedly committing the fatal sin of confusing numbers with their signs, when formally extending the system of positive natural numbers to include the negative ones. As mere symbols drawn on paper, without any content, the signs for negative numbers can have physical properties, but no arithmetical properties. Thus, Hankel's formal extension is really no arithmetical extension. But as we have seen above, Hankel argued, like Peacock, that all formal extensions of natural numbers must be capable of arithmetical interpretation, and he took the principle of permanence to guarantee such an interpretation, thereby showing that formal extensions are meaningful, rather than mere symbols, without any arithmetical content. Frege's criticism suggests that he doubted that Hankel's argument, along with his view of the principle of permanence, could be defended.

Frege seems to have thought in a similar manner about Peano's conception of permanence as a principle of practical rationality. It is well known that Peano criticized Frege's logical notation, introduced in the \textit{Begriffsschrift}, as deficient overall: from a scientific point of view, as he promptly noted in a review of the first volume of Frege's \textit{Grundgesetze}, ``the [\textit{Formulario}] system amounts to a more penetrating analysis. And then from the practical viewpoint, by using a composite sign to represent logical multiplication, Frege obscures its commutative and associative properties.'' (Peano 1895, 30) The criticism implies that a new notation must be such that these properties are not only preserved, but at least as perspicuous as in other already established notations. The ``composite'' signs for disjunction and conjunction, as constructed by Frege from the signs for implication and negation, triggered Peano's genuine concern that they could (as they in fact did) incur a practical cost, both with regard to the dissemination and the assimilation of Frege's logical calculus. Peano further contended that, unlike Frege's notation, the one adopted in his \textit{Formulario} is ``identical with those of Schr\"oder and Peirce'' (Peano 1895, 30), and added also that the deviation from Boole's notation is rather minimal and insignificant.\footnote{Schr\"oder disagreed and criticized Peano's notation as ``regrettably [...] a greatly diverging system of denotation'', and compared him to those who ``persist in still using sailing ships whilst steamboats have already been invented, constructed and are waiting at their service.'' (Schr\"oder 1898, 61) For historical details about Schr\"oder's general view on logical notations, see e.g. Peckhaus (1991). For discussion of Peano's notation, see Quine (1987). More recently, Frege's notation has been discussed in Belluci \textit{et al.} (2018) and Schlimm (2018).} One can discern here again Peano's conviction, which would later become more definitely articulated, that one should strive to make notations more uniform across science. In his responses, Frege emphatically rejected Peano's evaluation of the \textit{Begriffsschrift} from the scientific point of view, as less ``penetrating'' than that of the \textit{Formulario} (Frege 1896a, 1896b). Surprisingly, however, Frege chose not to respond at all to Peano's evaluation of the \textit{Begriffsschrift} from the \textit{practical} point of view, as ``obscure'', i.e., as making the properties of conjunction and disjunction less perspicuous than other systems of notation (Frege 1896b, 32). It's fair to say that Frege was rather unconcerned with practical matters such as increasing the uniformity of logical notations, and thought that the epistemological benefits of his notation would simply outweigh any possible practical costs with its assimilation and dissemination within the logic community. As is well-known, he was too optimistic in this regard.

Before Peano, Schr\"oder had strongly criticized the diagrammatic notation in the \textit{Begriffsschrift} as impractical, and more particularly, he dismissed what he perceived as Frege's alleged analogy between the diagrammatic notation and the usual arithmetical notation: ``it must be said that Frege's \textit{Begriffsschrift} promises too much in its title -- or more precisely, that the title does not correspond at all to the content. [...] In the subtitle [...] I find the very point in which the book corresponds least to its advertised program. [...]  If, to the impartial eye, the `modelling' appears to consist of nothing more than using \textit{letters} in both cases, then it seems to me this does not sufficiently justify the epithet used.'' (Schr\"oder 1880, 221) Thus, according to Schr\"oder, the fact that letters are used in both arithmetic and the \textit{Begriffsschrift} is not enough to make the former a notational model for the latter. 

Frege, of course, had meant the modeling in a different sense: ``That [the \textit{Begriffsschrift}] is modeled upon the formula language of arithmetic, as I indicated in the title, has to do with fundamental ideas rather than with details of execution (\textit{Einzelgestaltung}). [...] The most immediate point of contact between my formula language and that of arithmetic is the way in which the letters are employed.'' (Frege 1879, 6) So the fundamental idea is not that letters are employed, but \textit{how} they are employed. How, more exactly, are letters employed? Frege's view appears to be the following: ``I adopt the fundamental idea of distinguishing two kinds of signs [...] those whereby one can represent various things, and those that have a completely determined meaning. The former are the \textit{letters}, and these are to serve mainly for the expression of \textit{generality}.'' (Frege 1879, 10sq) This suggests that the kind of modeling that Frege had in mind is based on the use of letters as variables. But then, as Schr\"oder correctly pointed out, since the Boolean notation also uses letters for the expression of generality, Frege's criterion would not be enough to distinguish his logical notation.

The point that Schr\"oder missed here is, of course, that the fundamental ideas that Frege referred to are not related to the use of letters as variables, but to his reduction of arithmetical ordering to logical consequence, a reduction that he motivated in the following way: ``We divide all truths that require justification into two kinds, those for which the proof can be carried out purely by means of logic and those for which it must be supported by facts of experience. [...] Now, when I came to consider the question to which of these two kinds the judgments of arithmetic belong, I first had to ascertain how far one could proceed in arithmetic by means of inferences alone, with the sole support of those laws of thought that transcend all particulars. My initial step was to attempt to reduce the concept of ordering in a sequence to that of \textit{logical} consequence, so as to proceed from there to the concept of number. [...] This [...] led me to the idea of the present Begriffsschrift.'' (Frege 1879, 5sq) Thus, Frege's subtitle was justified by his attempted conceptual reduction of arithmetical ordering in a sequence to logical consequence. Hence, the modeling was meant as conceptual, rather than notational (Toader 2004). Notational modeling, motivated by rendering notations more uniform, did not concern Frege at all.

That Schr\"oder missed the point of Frege's conceptual modeling is quite clear from the way he concluded his review: ``The `appendix' of the \textit{Begriffsschrift} concerns `Some Topics from a General Theory of Sequences' and appears very abstruse -- the schemata are ornate with symbols! [...] The `sequence' is characterized only by the fact that a certain kind of advancement (which is otherwise left general) from one element to another is possible. [...] and the author is proud of the great generality that is given in this way to the concept of \textit{sequence}. It seems to me, however, that there is absolutely nothing of value in such a generalization.'' (Schr\"oder 1880, 230sq.)  It should be admitted, however, that beside their misunderstanding of Frege's foundational reasons that motivated his notation, Schr\"oder and Peano rejected that notation mainly because of a genuine concern with practical matters. In particular, and this is quite explicit in Peano's review, it was a commitment to permanence as a principle of practical rationality that justified skepticism with respect to a new logical notation that was so conspicuously dissimilar to other, more familiar notations.

\section{Beyond Peano}

Peano's argument, discussed above in section three, in favor of an explicitly Machian interpretation of permanence as a thought-economical principle of practical rationality, was reiterated by Hans Hahn in his criticism of the view according to which the principle of permanence is in fact logically necessary. The focus of this section will also be on his further point, that although arbitrariness cannot be entirely eliminated by insisting on the logical necessity of permanence, once the choice as to what formal laws to preserve and what to forgo when extending a mathematical system has been made, nothing else should remain arbitrary. Hahn discussed this in his review of a widely read book on analysis, \textit{Vorlesungen \"uber Zahlen- und Funktionenlehre}, published in 1916 by Alfred Pringsheim, who was described by Hahn as ``the most eminent representative of the arithmetical school, a school which goes so far in rejecting any geometrical element in analysis that it even renounces the help of the highly suggestive geometrical terminology which would make its propositions and proofs easier to understand.'' (Hahn 1919, 66)

Pringsheim had contributed to the same \textit{Encyklop\"adie} volume that included Schubert's article, discussed by Peano, a text titled \textit{Irrationalzahlen und Konvergenz unendlicher Prozesse}. This, too, attracted Frege's attention, who did not miss a chance to mock Pringsheim for claiming that ``the rational numbers feature as signs that \textit{may} represent well determined quantities but do not \textit{have} to.'' (Frege 1903, 83sq.) Frege saw this as yet another expression of the fatal sin he had seen committed by Hankel and Schubert (and many others): ``Clearly, this author too understands the rational numbers as the kind of figures that are artificially produced by a writing instrument on a writing surface or by a printing press.'' (Frege 1903, 83sq.) Pringsheim was reporting on a view he attributed to Weierstrass and Cantor -- ``a particular \textit{formal presentation} of irrational numbers,'' according to which an irrational number ``appears rather as a complete, \textit{newly created object}, or more concretely with Heine, as a \textit{new number sign} [...] with which one can calculate according to specific rules.''\footnote{``[...] eine bestimmte \textit{formale Darstellung} der Irrationalzahlen [...] erscheint vielmehr als ein fertiges, \textit{neu geschaffenes Objekt}, oder, noch konkreter nach Heine, als ein \textit{neues Zahlzeichen} [...] mit welchem nach bestimmten Regeln \textit{gerechnet} werden kann.'' (Pringsheim 1898, 54)} When mentioning Eduard Heine, Pringsheim duly quoted him as saying the following: ``I take a purely formal standpoint about the definition (of numbers) \textit{insofar as I call certain tangible signs numbers}, so that the existence of these numbers does not come into question then.''\footnote{``Ich stelle mich bei der Definition (der Zahlen) auf den rein formalen Standpunkt, \textit{indem ich gewisse greifbare Zeichen Zahlen nenne}, so dass die Existenz dieser Zahlen also nicht in Frage steht.'' (Pringsheim 1898, 54, n.21)} But Pringsheim immediately went on to clarify that these new number-signs should not be considered as standing for any quantities, and neither should the relations between them be seen as relations between quantities; rather, they should be seen as merely formal succession relations. He thought that this explained the sense in which a purely formal extension of the rational numbers to the irrational ones is to be understood: ``In particular, the concept of \textit{rational} numbers also undergoes an \textit{extension} in the \textit{sense} that they appear as signs to which initially only \textit{a certain succession} belongs, and which well \textit{can}, although do not {have to}, represent  particular \textit{quantities}.''\footnote{``Insbesondere erleidet hierbei also auch der Begriff der \textit{rationalen} Zahlen eine \textit{Erweiterung} in \textit{dem} Sinne, dass sie als Zeichen erscheinen, denen in erster Linie lediglich \textit{eine bestimmte Succession} zukommt, und die wohl bestimmte \textit{Quantit\"aten} vorstellen \textit{k\"onnen}, aber nicht \textit{m\"ussen}.'' (Pringsheim 1898, 55)} Unlike Hankel, Pringsheim thought that a purely formal extension was never at risk of being considered ``eine h\"ochst abstruse und beschwerliche K\"unstelei'', even if no arithmetical interpretation, whereby the new number-signs would come to stand for quantities, was available or even possible.\footnote{See footnote 4 on page 6 above for Pringsheim's criticism of Hankel on this point.} However, it seems fair to say that since he was merely  reporting on views about irrational numbers that had been propounded by others, it's not clear to what extent Pringsheim himself thought, at this time, that Heine's view was correct. 

In his review of Pringsheim's 1916 book, Hahn would also express doubts about the apparent identification of numbers and number-signs: ``Mr. Pringsheim certainly does not take the view that the figure 1 I am now drawing with ink on paper is a natural number: it is obvious that he does not mean anything so concrete but something else, though it would not be easy to give a precise formulation to the underlying thought.'' (Hahn 1919, 57) Of course, Hahn's point that what Pringsheim actually maintained was not clear enough is only slightly less damaging than Frege's point that he was simply wrong. But the real points of contention, as far as Hahn was concerned, were different: ``The few points on which we cannot completely follow the author concern individual questions of a methodological nature, of which we will only single out one for closer scrutiny, viz. the so-called principle of permanence.'' (Hahn 1919, 61). 
Before we can analyze these points, here is what Pringsheim wrote about permanence in his 1916 \textit{Vorlesungen}: ``To establish how these new numbers [i.e., zero, negative, rational, etc.] are to be \textit{ordered} or, alternatively, \textit{incorporated} into the ordered sequence of already existing numbers and how we are to \textit{calculate} with them, we shall make use of the \textit{transfer principle} which (following Hankel) is usually (but not very felicitously) called the principle of `\textit{permanence}', and we shall make use of it in what I regard as a notably improved form which bestows on it the character of a certain logical necessity. For in every case we shall introduce \textit{new number signs}, but only to such an extent that a subset of them represents signs for \textit{already existing} numbers. The latter are therefore already governed by certain rules establishing their succession and defining the arithmetical operations for them, and these rules can without further ado be \textit{transferred} into the new notation. If we are not to allow complete confusion in the manipulation of the total supply of newly created signs, we have hardly any choice but to extend the rules already governing part of it to the totality by definition, and to legitimize this step by proving that the stipulations we have made satisfy the requirements to be met by them without contradiction.''\footnote{Pringsheim (1916, VIIsq), quoted from the English translation of Hahn (1919, 61sq).} Following his contribution to the \textit{Encyklop\"adie}, Pringsheim considered the extension of a number system in a purely formal sense: he took the new number signs to have no interpretation, no reference to quantities, and the relations between them to be merely succession relations. Only some of the new signs are interpretable, since they can stand for numbers in the unextended system, but Pringsheim considered also this subset of the new signs as uninterpreted, and took the relations between them as formal succession relations. Based on this view of what it is to extend a number system, he argued that the formal rules that hold for the subset must be preserved or transferred to the entire set of new number signs, lest ``complete confusion'' be the outcome. This suggests that, according to him, if one dropped any rules, the formal succession relations between the new number signs could not be established, and the operations for manipulating these signs would fail to be properly defined. If this were the case, then calculating with the new numbers in the extended system would be logically impossible, which would defeat the very purpose of extension.

Against this view, Hahn recalled Peano's criticism of Schubert's attempt to formulate permanence precisely: ``Peano has demonstrated convincingly that this attempt is a failure: According to Schubert, the principle of permanence requires that the extensions of the number domain be carried out in such a way `that numbers in the   extended sense are governed by the \textit{same rules} as numbers in the non-extended sense'. But this requirement \textit{cannot possibly be met}. If it was really the case that numbers in the extended sense were still governed by the same rules, then they could not be distinguished from numbers in the non-extended sense, and we should not then be presented with an extension of the number domain. Hence it cannot be required that \textit{all} propositions of the original number domain continue to hold in the extended domain; the most that can be required is that the most important propositions continue to hold. But what the most important propositions are is a matter of personal judgement. The principle of permanence thereby ceases to be logical in nature and becomes at best a piece of methodological advice containing within itself an element of arbitrariness.'' (Hahn 1919, 62) Unlike Peano, Hahn appears to have misread Schubert's formulation of the principle of permanence: the third part of that formulation, which we have seen above, requires that the same propositions that hold for numbers in the extended sense are provably true for numbers in the non-extended sense. It does not require, as Hahn had it, that the same \textit{rules} that govern the numbers in the extended sense are the same \textit{rules} that govern numbers in the non-extended  sense. Although similar to Peano's argument, Hahn's version can be reconstructed as follows: (H1) It is not the case that all propositions that are true in the non-extended number system remain true in the extended one. (H2) Hence, only the most important propositions must be preserved. (H3) Thus, the principle of permanence is not a logically necessary principle. 

Whereas Peano had rejected the third part of Schubert's formulation of the principle of permanence, by denying that all propositions that are true in an extended system are provably true in the non-extended one, Hahn justified his (H1) as follows: if the rules that govern the extended system are the same as those that govern the non-extended one, then the two systems would be indistinguishable, which undercuts the claim that an extension is actually presented. Thus, Hahn inferred, it is not the case that all propositions that are true in the non-extended number system remain true in the extended one. But even if we take him to mean that the rules of the two systems are actually identical, Hahn's justification of (H1) does not fare any better than (and for the same reason noted when we discussed) Peano's justification of (P1) in his own argument. However, if one grants Hahn the justification for (H1), then qualifying the principle of permanence as suggested by his (H2) seems well motivated. The demand that only the most important propositions must be preserved was enough for Hahn to infer (H3), that the principle of permanence is not logically necessary, but is rather arbitrary.

However, Hahn did not follow Peano in thinking that permanence should, therefore, be reconceived as a thought-economical principle of practical rationality. He considered that ``there is more to be done than just insisting on the arbitrariness, that what is in question is not only an \textit{element of arbitrariness} but also an \textit{element of lawfulness} which must in turn be brought out into the open.'' (Hahn 1919, 63) These two elements, Hahn added ``cooperate in the usual extensions of the number domain [but] are not as sharply separated as they might be.'' (Hahn 1919, 63) Thus, in contrast to Peano, Hahn maintained that the source of arbitrariness should be determined more precisely. For he argued that even if permanence were a logically necessary principle, universally valid and applicable to all rules of the non-extended system, as Pringsheim took it to be, logical arbitrariness would still not be completely eliminated. Moreover, if one assumed that, when extending a system, once one decided what is the set of most important propositions to be preserved, all logical arbitrariness would subsequently disappear, and its place would be entirely taken by lawfulness, i.e., all other propositions in the extended system would follow logically from that set, one would simply be making a false assumption. 

Here is Hahn's argument. Following Pringsheim, let us consider the formal extension of the natural numbers by positive rational numbers. Introduce a set of entirely new signs for positive rational numbers, and let a subset of them be new signs for natural numbers. Let us distinguish ``proper'' fractions (whose numerators are not multiples of their denominators) from ``improper'' fractions (whose numerators are multiples of their denominators). In order to proceed with the extension, one must transfer the formal rules for improper fractions (which are just the ``direct consequences'' of the rules for natural numbers) to the proper fractions, thereby establishing the formal succession relations between proper fractions and defining the basic operations for them. This allows us to take proper fractions as signs for positive rational numbers, and it allows us to calculate with them. But this can be achieved, Hahn argued, only if the formal rules for improper fractions can be extended so as to define equality between proper fractions in the same way as equality between improper ones, that is, 

\begin{center}
$\frac{b}{a} = \frac{b'}{a'}$ iff $ba' - ab' = 0$.    
\end{center}

If the extension of the natural numbers proceeded along these lines, and if equality could be uniquely defined as above, then ``all arbitrariness would be excluded'' (Hahn 1919, 64). But Hahn denied that Pringsheim’s account of extension can guarantee the complete elimination of arbitrariness. For he correctly maintained that equality between improper fractions is not uniquely defined. As he put it, ``we have just as much right to express'' equality between improper fractions, for example, in the following way:

\begin{center}
$\frac{b}{a} = \frac{b'}{a'}$ iff $(ba' - ab')^{2} + (r_{a}(b) - r_{a'}(b'))^{2} = 0$,
\end{center}

where $r_{a}(b)$ and $r_{a'}(b')$ are the absolutely smallest residues of \textit{b} and \textit{b'} to moduli \textit{a} and \textit{a'}, respectively. For improper fractions, these are equivalent definitions of equality because $(r_{a}(b) - r_{a'}(b'))$ is always zero, but when transferred to proper fractions this may not always be zero, so the two definitions are not equivalent any longer. This entails, according to Hahn, that there is more than one logically possible way of extending the natural numbers to positive rationals: ``If it were true that the domain of numbers could only be extended by assuming that [the first definition of equality] continued to hold, we should have just as much right to make the same assumption about [the second definition of equality]. But this would lead to an entirely different definition of equality for improper fractions. Therefore it is evidently not the case that we may conclude without further ado that [the first definition of equality] gives us the only possible definition of equality for improper fractions. And if we prefer the definition of equality given us by [the first definition of equality] to the one given us by [the second definition of equality], which in and by itself is just as possible, it is not, it seems, because we are under a logical compulsion. Again, arbitrariness seems to reign supreme.'' (Hahn 1919, 64) 

While he further noted that practical considerations like ``realizability and applicability'' (Hahn 1919, 71, n.14) may sometimes justify our preferences as to one logically possible extension of a number system over another, Hahn argued that, in the case of positive rational numbers, uniqueness of extension can be recovered if one stipulates that the basic operations must be one-to-one [\textit{eindeutig umkehrbar}] mappings (\textit{Verkn\"upfungen}). Hence, the task as he saw it: ``the system of natural numbers is to be extended by adding new `numbers' in such a way that multiplication [and thus division] will always be one-to-one in the extended system.'' (Hahn 1919, 63) 

Here is how that stipulation works. Consider the extended system of positive rational numbers, including division and multiplication as basic operations. In order to show, for any proper fractions $\frac{b}{a}$ and $\frac{b'}{a'}$, that their equality can be uniquely defined, assume $\frac{b}{a} = \frac{b'}{a'}$. Since multiplication is assumed one-to-one, we then have $\frac{b}{a}(aa') = \frac{b'}{a'}(aa')$, but because multiplication is associative and commutative, we have $(\frac{b}{a}a)a' = (\frac{b'}{a'}a')a$, and thus $ba' = ab'$. Conversely, assume that $ba' = ab'$. Since division is one-to-one, we have $ba'(\frac{a}{a}) = ab'(\frac{a'}{a'})$, and because division is associative and commutative, $\frac{b}{a}(aa') = \frac{b'}{a'}(aa')$ and finally, $\frac{b}{a} = \frac{b'}{a'}$. Thus, if the basic operations for proper fractions are taken to be one-to-one mappings, then uniqueness is restored and logical arbitrariness is fully eliminated by means of exclusively logical (rather than practical) considerations.

\section{Conclusion}

I have focused on Peano's argument for his conception of the principle of permanence as a principle of practical rationality --- a Machian conception, driven by the ideal of thought economy, that Peano believed motivated Burali-Forti and Marcolongo's notational project for vectorial calculus, and that he also summoned against Frege's diagrammatic notation for logic. Unlike earlier mathematicians, such as Peacock and Hankel, who considered the principle of permanence a principle of theoretical rationality, since it was meant to guarantee the arithmetical interpretability and thus the meaningfulness of symbolic algebra and formal mathematics, Peano thought that its less than universal validity made it an arbitrary principle, which he further thought implied that it could only be understood as a principle of practical rationality, especially relevant in the selection of more familiar and uniform notations. Insistence on its logical necessity, by mathematicians like Schubert and Pringsheim, led Hahn to revive Peano's argument, and to identify a further source of arbitrariness, the elimination of which, as Hahn pointed out, requires further logical considerations.

\section{References}

Bellucci, F., A. Moktefi and A.-V. Pietarinen, 2018. Simplex sigillum veri: Peano, Frege, and Peirce on the Primitives of Logic. History and Philosophy of Logic, 39, 80--95. 

Burali-Forti, C. and R. Marcolongo, 1909. Elementi di Calcolo Vettoriale. Con numerose applicazioni alla geometria, alla meccanica, e alla fisica matematica. 
Seconda edizione  riordinata e ampliata, Nicola Zanichelli, Bologna, 1921.

De Morgan, A., 1836-1842. The Differential and Integral Calculus. Baldwin and Cradock.

De Morgan, A., 1837. Elements of Algebra, Preliminary to the Differential Calculus. 2nd ed., Taylor and Walton.

Detlefsen, M., 2005. Formalism. In: Shapiro, Stewart (Ed.), The Oxford Handbook of Philosophy of Mathematics and Logic. Oxford University Press, 236--317.

Fisch, M., 1999. The Making of Peacock's Treatise on Algebra: A Case of Creative Indecision. Archive for History of Exact Sciences, 54, 137--179.

Frege, G., 1879. Begriffsschrift: eine der arithmetischen nachgebildete Formelsprache des reinen Denkens. Halle, Louis Nebert. Translated as Begriffsschrift, a formula language, modeled upon that of arithmetic, of pure thought. In van Heijenoort, J. (Ed.), From Frege to Gödel. A Source Book in Mathematical Logic, 1879-1931. 4th ed., Harvard University Press, 1--82.

Frege, G., 1896a. \"Uber die Begriffsschrift des Herrn Peano und meine eigene. Berichte \"uber die Verhandlungen der K\"oniglich S\"achsischen Gesellschaft der Wissenschaften zu Leipzig. Mathematisch-Physische Klasse, 48, 361--378. Translated in Collected Papers in Mathematics, Logic, and Philosophy. Basic Blackwell, 1984, 234--248.

Frege, G., 1896b. Lettera del sig. G. Frege all'Editore. Rivista di Matematica, 6, 53--61. Translated in Southern Journal of Philosophy, 9, 32--36. 

Frege, G., 1899. \"Uber die Zahlen des Herrn H. Schubert, H. Pohle, Jena. Eng. tr. in Collected Papers in Mathematics, Logic, and Philosophy. Basic Blackwell, 1984, 249--272.

Frege, G., 1903. Die Grundgesetze der Arithmetik, II, Jena. Translated as Basic Laws of Arithmetic, I and II, Oxford University Press, 2013.

Hahn, H., 1919. Review of Alfred Pringsheim: Vorlesungen \"uber Zahlen- und Funktionenlehre. Göttingische Gelehrte Anzeigen, 9, 321--347. Translated in Empiricism, Logic, and Mathematics. Philosophical Papers. D. Reidel, 1980, 51--72.

Hankel, H., 1867. Theorie der complexen Zahlensysteme. Leopold Voss, Leipzig.

Hollings, C., U. Martin, and A. Rice, 2017. The Lovelace–De Morgan mathematical correspondence: A critical re-appraisal. Historia Mathematica, 44, 202--231.

Koppelman, E., 1971. The Calculus of Operations and the Rise of Abstract Algebra. Archive for History of Exact Sciences, 8, 155--242.

Mach, E., 1905. Erkenntnis und Irrtum. J. A. Barth. Translated as Knowledge and Error. D. Reidel, 1976.

Musil, R., 1908. Beitrag zur Beurteilung der Lehren Machs. PhD dissertation at Berlin University. Translated as On Mach's Theories. Philosophia Verlag, Munich, and The Catholic University of America Press, Washington, 1982.

Peacock, G., 1833. Report on the Recent Progress and Present State of Certain Branches of Analysis. Report of the British Association, 185--352.

Peacock, G., 1842-1845. A Treatise on Algebra. Vol. 1 and 2, London.

Peano, G., 1895. Review of G. Frege, Grundgesetze der Arithmetik, begriffs
schriftlich abgeleitet, vol. I. Rivista di Matematica, 5, 122--128. Translated in Southern Journal of Philosophy, 9, 27--31.

Peano, G., 1903. Principio de Permanentia. Exercitio de Latino recto. Rivista di matematica, VIII, 84--87.

Peano, G., 1910. Sui fondamenti dell' analisi. Mathesis Societa Italiana di Matematica Bolletino, 2, 31--37. Translated in Peano 1973, 219--226.

Peano, G., 1911. Le definizioni in matematica. Arxivs de L'Institut de Ciencies, 1, 49--70. 

Peano, G., 1921. Le definizioni in matematica. Periodico di Mathematiche, 4, 175--189. Translated in Peano 1973, 235--246.

Peano, G., 1973. Selected Works of Giuseppe Peano, George Allen \& Unwin, London.

Peckhaus, V., 1991. Ernst Schröder und die `pasigraphischen Systeme' von Peano und Peirce. Modern Logic, 1, 174--205.

Peckhaus, V., 1997. Logik, Mathesis universalis und allgemeine Wissenschaft. Leibniz und die Wiederentdeckung der formalen Logik im 19. Jahrhundert. Akademie Verlag.

Peirce, C. S., 1877. The fixation of belief. In: Buchler, J (Ed.),  Philosophical Writings of Peirce. Dover, 1955, 5--21.

Peirce, C. S., 1878. How to make our ideas clear. In: Buchler, J (Ed.),  Philosophical Writings of Peirce. Dover, 1955, 23--41.

Pringsheim, A., 1898. Irrationalzahlen und Konvergenz unendlicher Prozesse. Encyklop\"adie der mathematischen Wissenschaften mit Einschluss ihrer Anwendungen, I. Teubner, Leipzig, 47--146. 

Pringsheim, A., 1916. Vorlesungen \"uber Zahlen- und Funktionenlehre, I. Teubner, Leipzig.

Pycior, H. M., 1981. George Peacock and the British Origins of Symbolic Algebra. Historia Mathematica, 8, 23--45.

Quine, W. V. O., 1987. Peano as Logician. History and Philosophy of Logic, 8, 15--24.

Sallent Del Colombo, E., 2010. Il dibattito sull’unificazione delle notazioni vettoriali. Il contributo di Cesare Burali-Forti e Roberto Marcolongo. In Roero, C. S. (Ed.), Peano e la sua scuola. Collana di Studi e Fonti, Centro per la storia dell'Università di Torino, 509--529.

Schlimm, D., 2018. On Frege's \textit{Begriffsschrift} Notation for Propositional Logic: Design Principles and Trade-Offs. History and Philosophy of Logic, 39, 53--79.

Schr\"oder, E., 1880. Review of Frege's \textit{Begriffsschrift}.  Zeitschrift f\"ur Mathematik und Physik, 25, 81--94. Translated in Conceptual Notation and related articles. Clarendon Press, Oxford, 1972, 218--232.

Schr\"oder, E., 1898. On Pasigraphy. Its Present State And The Pasigraphic Movement In Italy. The Monist, 9, 44--62.

Schubert, H., 1898. Grundlagen der Arithmetik. Encyklop\"adie der mathematischen Wissenschaften mit Einschluss ihrer Anwendungen. Teubner, Leipzig, 1--27.

Schubert, H., 1899. Elementare Arithmetik und Algebra. Göschen, Leipzig.

Schubring, G., 2005. Conflicts between Generalization, Rigor, and Intuition. Number Concepts Underlying the Development of Analysis in 17--19th Century France and Germany. Springer.

Tappenden, J., 2019. Infinitesimals, Magnitudes, and Definition in Frege. In: Ebert, P. A. and M. Rossberg (Eds.), Essays on Frege's Basic Laws of Arithmetic. Oxford University Press, 235--263. 

Toader, I. D., 2004. On Frege's Logical Diagrams. Diagrammatic Representation and Inference. Lecture Notes in Computer Science, 2980, Springer, 22--25.

Toader, I. D., 2019. Talking Past Each Other: Mach and Husserl on Thought Economy. In: Stadler, Friedrich (Ed.), Ernst Mach – Life, Work, Influence. Springer, 213--221.

Wilson, E. B., 1913. The Unification of Vectorial Notations. Bulletin of the American Mathematical Society, 19/10, 524--530.

\end{document}